\documentclass{amsart}
\usepackage{amsmath,amssymb,amscd,amsthm}
\usepackage{tikz}
\usetikzlibrary{positioning}
\usepackage[pagebackref,colorlinks=false,urlbordercolor=white]{hyperref}

\newcommand{\ov}{\overline}
\newcommand{\cM}{\mathcal M}
\newcommand{\rsp}{\raisebox{0em}[2.7ex][1.3ex]{\rule{0em}{2ex} }}
\newtheorem{lem}{Lemma}
\newtheorem{prop}[lem]{Proposition}
\newtheorem{thm}[lem]{Theorem}

\newtheorem{cor}[lem]{Corollary}

\newcommand{\Z}{{\mathbb Z}}

\newcommand{\Q}{{\mathbb Q}}

\newcommand{\hra}{\hookrightarrow}
\newcommand{\lra}{\longrightarrow}

\newcommand{\Sel}{{\operatorname{Sel}}}
\newcommand{\Gal}{{\operatorname{Gal}}}
\newcommand{\Ver}{{\operatorname{Ver}}}
\newcommand{\MD}{{\operatorname{M}}}

\newcommand{\Aut}{{\operatorname{Aut}}}

\newcommand{\SD}{{\operatorname{SD}}}
\newcommand{\Hol}{{\operatorname{Hol}}}
\newcommand{\Cl}{{\operatorname{Cl}}}
\newcommand{\Am}{{\operatorname{Am}}}
\newcommand{\gen}{{\operatorname{gen}}}
\newcommand{\disc}{{\operatorname{disc}\ }}

\newcommand{\la}{\langle}
\newcommand{\ra}{\rangle}
\newcommand{\eps}{\varepsilon}
\newcommand{\vr}{\varrho}

\newcommand{\fa}{\mathfrak a}

\newcommand{\fp}{\mathfrak p}

\newcommand{\cO}{\mathcal O}

\title[Extensions with Galois group $\Hol(C_8)$]
      {Extensions with Galois group $\Hol(C_8)$ unramified over a
       complex quadratic number field}

\author{E. Benjamin, F. Lemmermeyer, C. Snyder}

\begin{document}

\maketitle

\begin{abstract}
  Nous \'etudions des extensions normales $L/\Q$ avec groupe de Galois
  $\Hol(C_8)$ qui sont non-ramifi\'es sur un sous-corps quadratique
  complexe $K$. Le groupe de Galois $\Gal(L/K)$ est soit le groupe
  quasidi\'edral ou le groupe modulaire d'ordre $16$. Nous pr\'esentons
  une construction explicite de tels corps.

  \medskip
  
  We study normal extensions $L/\Q$ with Galois group $\Hol(C_8)$
  that are unramified over a complex quadratic subfield $K$. The
  Galois group $\Gal(L/K)$ is either the semi-dihedral group
  or the modular group of order $16$. We present
  an explicit construction of such fields.
\end{abstract}

Let $L$ be a normal extension of $\Q$ with Galois group $\Hol(C_8)$,
the holomorph of the cyclic group $C_8$ of order $8$, i.e.,
\begin{align*}
  \Hol(C_8) & \simeq C_8 \rtimes \Aut(C_8) \\
            & = \la a, x, y: a^8 = x^2 = y^2 = (xy)^2 = 1, xax^{-1} = a^{-1},
    yay^{-1} = a^5 \ra. 
\end{align*}
Assume that such an extension $L$ is unramified over some quadratic
subextension $k/\Q$ of $L/\Q$. We will show that
$\Gal(L/k) \simeq \SD_{16}$ or $\Gal(L/k) \simeq \MD_4(2)$. Here $\SD_{16}$
is the semi-dihedral group of order $16$, and $\MD_4(2)$ is the modular
group of order $16$. These groups have the following presentations:
\begin{align*}
  \SD_{16} & = \la \sigma, \tau: \sigma^8 = \tau^2 = 1, \
                 \tau \sigma\tau = \sigma^3 \ra, \\
  \MD_4(2) & = \la \sigma, \tau: \sigma^8 = \tau^2 = 1,  
                \tau \sigma \tau = \sigma^5 \ra. 
\end{align*}

In this article we will classify all complex quadratic number fields
whose $2$-class field tower terminates at the second step $k^2$ and
for which $\Gal(k^2/k)$ is either $\SD_{16}$ or $\MD_4(2)$, and we will
show how to construct these extensions explicitly. The case of real
quadratic fields requires more effort and will be taken care of later.

\section{Group Theoretical Calculations}

Let us consider the following situation. Let $F$ be a number field
with odd class number, $L/F$ a normal $2$-extension with Galois group
$G = \Gal(L/F)$, and $k/F$ a quadratic subextension; then $H = \Gal(L/k)$
is a normal subgroup of $G$ with index $2$. By Chebotarev's monodromy
theorem, $G$ is generated by the inertia subgroups $T_\fp(L/F)$ of prime
ideals $\fp$ in $F$; thus if $L/k$
is unramified, then $G$ is generated by elements of order $2$. In fact in
this case we have $1 = T_\fp(L/F) \cap H$, hence $G$ is generated by
elements of order $2$ lying outside of $H$.

\begin{prop}\label{Pr1}
Let $F$ be a number field with odd class number, $L/F$ a normal
$2$-extension with Galois group $G$, $H$ a subgroup of index $2$, and
let $k = L^H$ denote the fixed field of $H$. If  $L/k$ is unramified, then
$G$ is generated by elements of order $2$ lying outside of $H$.
\end{prop}

Calculations with GAP show:

\begin{thm}\label{Th1}
  The following groups $G$ are all the nonabelian groups of order $16$
  that possess subgroups $H$ of index $2$ such that $G$ is generated by
  the elements of order $2$ in $G \setminus H$:  
  $$ \begin{array}{ccc|ll}
    G & \text{HS} & & H &   \\ \hline
    \rsp 16.07 & 16.012 & D_8                & 8.01 &  C_8            \\
    \rsp 16.11 & 16.006 & C_2 \times D_4     & 8.02 &  C_4 \times C_2 \\
    \rsp 16.11 & 16.006 & C_2 \times D_4     & 8.03 &  D_4            \\
    \rsp 16.13 & 16.008 & C_4 \curlyvee D_4  & 8.04 &  H_8           
  \end{array} $$
  Here the groups are identified by their GAP ID as well as by their
  number in the tables of Hall \& Senior.
  
  Thus if $L/F$ is a normal nonabelian extension of degree $16$ of a
  field $F$ with odd class number, and if $L$ is unramified over some
  quadratic subextension $k/F$, then $G = \Gal(L/F)$ is one of the
  four groups in this table.
\end{thm}

For nonabelian groups $G$ of order $32$, there are, of course, a lot more
possibilities. A similar calculation with GAP yields:

\begin{thm}
  The following groups $G$ are all the nonabelian groups of order $32$ that
  possess subgroups $H$ of index $2$ such that $G$ is generated by the
  elements of order $2$ in $G \setminus H$:
  $$ \begin{array}{ccc|ll}
    G & & \text{HS} & H &  \\ \hline
   \rsp 32.18 & D_{16} & 32.049 & 16.01 &  C_{16}    \\
   \rsp 32.27 & & 32.033 & 16.03 &  (C_4 \times C_2) \rtimes C_2  \\
   \rsp 32.28 & & 32.036 & 16.04 &  C_4 \rtimes C_4  \\
   \rsp 32.34 & D(4,4) & 32.034 & 16.02 &  C_4 \times C_4  \\
   \rsp 32.39 & C_2 \times D_8 & 32.023 & 16.07 &  D_{8}  \\
   \rsp    & &        & 16.05 &  C_8 \times C_2 \\
   \rsp 32.42 & & 32.026 & 16.09 &  H_{16}  \\   
   \rsp 32.43 & \Hol(C_8) & 32.044 & 16.08 &  \SD_{16}  \\   
   \rsp    & &        & 16.06 &  \MD_4(2)  \\
   \rsp 32.46 & D_4 \times V_4 & 32.008 & 16.11 &  C_2 \times D_4   \\
   \rsp    & &        & 16.10 &  C_4 \times C_2 \times C_2   \\
   \rsp 32.48 & C_2 \times (D_4 \curlyvee C_4)
            & 32.010 & 16.13 &  (C_4 \times C_2)  \rtimes C_2  \\
   \rsp    & &        & 16.12 &  C_2 \times H_8  \\
   \rsp 32.49 & & 32.042 & 16.11 &  C_2 \times D_4  \\
   \rsp    & &        & 16.13 &  (C_4 \times C_2) \rtimes C_2 \\
   \rsp 32.50 & & 32.043 & 16.12 &  C_2 \times H_8  \\
  \end{array} $$
\end{thm}

The groups $G$ and $H$ in these two tables are identified by their GAP ID; the
third column (HS) gives the number of $G$ in the tables by Hall and Senior.
The case we are interested in is the group $\Hol(C_8)$, for which the
table above yields the following result:

\begin{prop}
  Let $F$ be a number field with odd class number. If $L/F$ is a normal
  extension with Galois group $\Hol(C_8)$, and if $L/k$ is unramified for
  a quadratic subextension $k/F$ of $L/F$, then $\Gal(L/k) \simeq \SD_{16}$
  or $\Gal(L/k) \simeq \MD_4(2)$.
\end{prop}

In this article we will show that both cases occur, and we will show how to
explicitly construct such extensions in the special case where $F = \Q$ and
$k$ is a complex quadratic number field. Since in our case these extensions
$L/k$ coincide with the second Hilbert $2$-class field of $k$, they
are normal over $\Q$, and this implies that $\Gal(L/\Q) \simeq \Hol(C_8)$ in
both cases.

\subsection*{Determination of $\Gal(k^2/\Q)$}
Here we will determine the Galois groups $\Gal(k^2/\Q)$ for
unramified extensions $k^2/k$ of quadratic number fields $k$ with
$\Gal(k^2/k) \simeq \SD_{2^n}$ and $\Gal(k^2/k) \simeq M_n(2)$.

We start by looking at the classification problem group theoretically.
Both $\SD_{2^n}$ and $\MD_n(2)$ are semidirect products $C_{2^{n-1}}\rtimes C_2$,
where $C_2=\la \tau\ra$ acts on $C_{2^{n-1}}=\la \sigma\ra$ by
$$\tau\sigma\tau=\left\{
                   \begin{array}{ll}
                     \sigma^{2^{n-2}-1}, & \hbox{in}\;\;\SD_{2^n}, \\
                     \sigma^{2^{n-2}+1}, & \hbox{in}\;\;\MD_n(2).
                   \end{array}
                 \right. $$

Thus, letting $m=2^{n-2}$, we have the following presentations:
$$\SD_{2^n}=\SD_{4m}=\la \sigma,\tau: \sigma^{2m}=\tau^2=1,
   \tau\sigma\tau=\sigma^{m-1}\ra,$$
and
$$\MD_n(2)=\la \sigma,\tau: \sigma^{2m}=\tau^2=1,
  \tau\sigma\tau=\sigma^{m+1}\ra.$$

Now, if $\Gal(k^2/k)\simeq \SD_{2^n}$ or $\MD_n(2)$, then
$\Gal(k^2/\Q)$ is uniquely determined; in light of Proposition~1 this
follows from the following group theoretic result:

\begin{thm}\label{DT1}
  Let $H=\SD_{2^n}$ or $H=\MD_n(2)$. Then (up to isomorphism) there
  exists a unique finite $2$-group $G$ such that $H$ is a maximal
  subgroup of $G$ and $G$ is generated by elements of order $2$ not in
  $H$. Namely, letting $m=2^{n-2}$,
  $$G\simeq C_{2^{n-1}}\rtimes \Aut(C_{2^{n-1}})[2]$$
  $$=\la a,x,y:a^{2m}=x^2=y^2=(xy)^2=1,xax=a^{-1},yay=a^{m+1}\ra\,.$$
  Hence for both $H=\SD_{2^n}$ and $H=\MD_n(2)$, the groups $G$ coincide.
\end{thm}

\begin{proof}
  Existence: Let $G\simeq C_{2^{n-1}}\rtimes \Aut(C_{2^{n-1}})[2]$ as
  presented above. Then notice that $H=\la a,xy\ra\simeq \SD_{2^n}$
  and that $G=\la x,y, xa\ra$ is generated by elements of order $2$
  not in $H$.  Moreover $H=\la a, y\ra \simeq \MD_n(2)$ and $G=\la
  x,xy,xa\ra$.

  Uniqueness: Consider $H=\SD_{2^n}=\SD_{4m}=\la \sigma,\tau:
  \sigma^{2m}=\tau^2=1,\tau\sigma\tau=\sigma^{m-1}\ra.$ Suppose $G$ is
  a finite $2$-group with maximal subgroup $H$, say $G=H\cup\varrho H$
  for any $\varrho\in G\setminus H$, and such that
  $G = \la \varrho \eta : \eta\in H, |\varrho\eta| = 2\ra.$ Here
  $|\cdot|$ denotes the order of an element of a finite group.
  Without loss of generality,   pick $\varrho\in G\setminus H$ such
  that $\varrho^2=1.$

We will show that $\varrho\sigma\varrho=\sigma^{-1}$ and that,
moreover, $\vr$ and $\tau$ may be chosen so that
$\varrho\tau\varrho=\tau$, which will establish our uniqueness result
for this $H$. Toward this end, consider $\varrho\sigma\varrho$ and
notice that $|\varrho\sigma\varrho|=|\sigma|=2m$. Now, recall that for
integers $i$ and $j$, $|\sigma^i|=2m/(i,2m)$ and
$|\sigma^j\tau|=4/(j,2)$ (for notice that
$(\sigma^j\tau)^2=\sigma^j\tau\sigma^j\tau=\sigma^j\sigma^{j(m-1)}=\sigma^{mj}$).
The only elements in $H$ of order $2m$ are those of the form $\sigma^i$
for which $i$ is odd. Hence $\varrho\sigma\varrho=\sigma^u$ for some
odd integer $u$. Moreover, since $\sigma=\varrho^2\sigma\varrho^2=\sigma^{u^2}$,
we have $u^2\equiv 1\bmod 2m$.

Next, consider $\varrho\tau\varrho$ which has order $2$. From
above, $\varrho\tau\varrho=\sigma^m$ or $\sigma^{2v}\tau$ for some
integer $v$. But $\tau=\varrho^2\tau\varrho^2$ implies that
$\varrho\tau\varrho=\sigma^{2v}\tau$ and therefore
$\tau=\varrho^2\tau\varrho^2=\varrho\sigma^{2v}\tau\varrho=\sigma^{2v(u+1)}\tau$,
so that $2v(u+1)\equiv 0\bmod 2m.$

Thus so far we have determined that there are integers $u$ and $v$ such that
$$\vr\sigma\vr=\sigma^u\;(u^2\equiv 1\bmod
2m)\;\;\text{and}\;\;\vr\tau\vr=\sigma^{2v}\tau \;(2v(u+1)\equiv
0\bmod 2m).$$ Given these actions, we now need to determine those
$\eta\in H$ for which $|\vr\eta|=2$. For any $\eta\in H$, either
$\eta=\sigma^i$ or $\eta=\sigma^j\tau$ for some integers $i,j$. But
then $|\vr\sigma^i|=2$ iff $1=\vr\sigma^i\vr\sigma^i=\sigma^{i(u+1)}$
iff $i(u+1)\equiv 0\bmod 2m$. On the other hand, $|\vr
\sigma^j\tau|=2$ iff
$1=\vr\sigma^j\tau\vr\sigma^j\tau=\sigma^{j(u-1+m)+2v}$ iff
$j(u-1+m)+2v\equiv 0\bmod 2m.$ Thus we get
$$\la \vr\eta: \eta\in H, |\vr\eta|=2\ra=\la \vr, \sigma^i, \sigma^j\tau:i(u+1),j(u-1+m)+2v\equiv 0\bmod 2m\ra.$$
But now this group must be $G$. This is only possible if there is an
$i$ which is odd. Hence $u+1\equiv 0\bmod 2m$. This implies that
$\vr\sigma\vr=\sigma^{-1}$ as predicted above. Moreover, notice that
with $u\equiv -1\bmod 2m$, there are integers $j$ that satisfy
$j(m-2)+2v\equiv 0\bmod 2m$, since $m-2\equiv 2\bmod 4.$ This implies
that there's an integer such that $\sigma^j\tau \in G$ and thus
$\tau\in G$, too. Hence we have
$$ G=\la \vr,\sigma,\tau:\vr^2=\tau^2=\sigma^{2m}=1,
\vr\sigma\vr=\sigma^{-1},\vr\tau\vr=\sigma^{2v}\tau,
\tau\sigma\tau=\sigma^{m-1}\ra.$$
If $v$ is even, then let $\tau'=\sigma^v\tau$  (and then rename
$\tau'$ as $\tau$) and we obtain
$$ G = \la \vr,\sigma,\tau:\vr^2=\tau^2=\sigma^{2m}=1,
   \vr\sigma\vr=\sigma^{-1},\vr\tau\vr = \tau,\tau\sigma\tau=\sigma^{m-1}
   \ra. $$
Similarly, if $v$ is odd, then letting $\vr'=\vr\sigma^{-v+m/2}$ and
$\tau'=\sigma^{2v}\tau$, yields this $G$ again.
Thus $G\simeq C_{2^{n-1}}\rtimes \Aut(C_{2^{n-1}})[2]$, as is easily seen.

\bigskip

Now consider $H=\MD_{n}(2)=\la \sigma,\tau:
\sigma^{2m}=\tau^2=1,\tau\sigma\tau=\sigma^{m+1}\ra$ (where again
$m=2^{n-2}\geq 4$). Suppose $G$ is a finite $2$-group with maximal
subgroup $H$, say $G=H\cup\varrho H$ for any $\varrho\in G\setminus
H$, and such that $G=\la \varrho \eta:\eta\in H, |\varrho\eta|=2\ra.$
Without loss of generality, pick $\varrho\in G\setminus H$ such that
$\varrho^2=1.$

As in the previous case, we wish to show that
$\vr\sigma\vr=\sigma^{-1}$ and that in addition $\vr$ and $\tau$ may
be chosen so that $\vr\tau\vr=\tau$. Since $\sigma$ has order $2m$ and
$\tau$ has order $2$, we'll start by computing the orders of elements
$\eta$ in $H$.  As before, we have $|\sigma^i|=2m/(i,2m)$. On the
other hand, we can easily see that
$(\sigma^j\tau)^i=\sigma^{j(m+i)}\tau^i$ and thus in particular
$(\sigma^j\tau)^2=\sigma^{j(m+2)}$ which has order $m/(j,m)$ and thus
$|\sigma^j\tau|=2m/(j,m)$. We therefore have the following results for
$\eta\in H$:
\begin{align*} |\eta|=2m\;&\Leftrightarrow\;
  (\eta=\sigma^i\;\text{or}\; \eta=\sigma^j\tau\;
  \text{for some odd integers}\;i,j),\\
  |\eta|=2 \;&\Leftrightarrow\; (\eta=\sigma^m\;\text{or}\;
  \eta=\sigma^{\ell m}\tau\;\text{for some}\;\ell=0\;\text{or}\;1).
\end{align*}

Hence since $|\vr\tau\vr|=2$, we must have $\vr\tau\vr=\sigma^m$ or
$\vr\tau\vr=\sigma^{\ell m}\tau$. Moreover $\vr\sigma\vr=\sigma^i$ or
$\vr\sigma\vr=\sigma^j\tau$ for some odd $i,j$. But we claim that
$\vr\tau\vr=\sigma^m$ cannot happen. For otherwise,
$\tau=\vr^2\tau\vr^2=\vr\sigma^m\vr\in\la \sigma^m\ra$, as can be seen
by checking the two possibilities for $\vr\sigma\vr$.

Therefore, $\vr\tau\vr=\sigma^{\ell m}\tau.$ Moreover,
$\vr\sigma\vr=\sigma^u$ or$=\sigma^v\tau$ for some odd $u,v$. We now
claim that $\vr\sigma\vr=\sigma^v\tau$ cannot happen. For assume
$\vr\sigma\vr=\sigma^v\tau$ with $v$ odd. We will show that
$\sigma\not\in\la\vr\eta:\eta\in H,|\vr\eta|=2\ra.$ By some easy
calculations, we can see that
$$|\vr\sigma^i|=2\;\Leftrightarrow\;
  (i(v+1)\equiv m\bmod 2m,\;i\;\text{even}),$$
$$|\vr \sigma^j\tau|=2\;\Leftrightarrow\;
   (j(v+1)\equiv m(\ell+1)\bmod 2m,\;j\;\text{even}).$$
But then we see $\sigma\not\in\la\vr\eta:\eta\in H,|\vr\eta|=2\ra.$

\medskip
Hence we must have $\vr\tau\vr=\sigma^{\ell m}\tau$ and
$\vr\sigma\vr=\sigma^u$ for some odd integer $u$.  Since
$\sigma=\vr^2\sigma\vr^2$, a straightforward calculation shows that
$u^2\equiv 1\bmod 2m$. Then we have
$$|\vr\sigma^i|=2\;\Leftrightarrow\; i(u+1)\equiv 0\bmod 2m,$$
$$|\vr \sigma^j\tau|=2\;\Leftrightarrow\; j(u+1+m)\equiv \ell m\bmod 2m.$$
In order for the elements of order $2$ in $G\setminus H$ to generate $G$,
at least one of the integers $i$ must be odd. Hence $u\equiv -1\bmod 2m$,
in which case $j$ can be any  integer with $j\equiv \ell\bmod 2$. Thus we have
$$G = \la \vr,\sigma,\tau:\vr^2=\tau^2=\sigma^{2m}=1,
   \vr\sigma\vr=\sigma^{-1},\vr\tau\vr=
    \sigma^{\ell m}\tau,\tau\sigma\tau=\sigma^{m+1}\ra.$$
If $\ell=1$, then let $\vr'=\vr\sigma$ (and then replace $\vr'$ by $\vr$).
Hence in either case
$$G=\la \vr,\sigma,\tau:\vr^2=\tau^2=\sigma^{2m}=1,
   \vr\sigma\vr=\sigma^{-1},\vr\tau\vr=
    \tau,\tau\sigma\tau=\sigma^{m+1}\ra. $$
This proves the theorem.
\end{proof}

\begin{cor}\label{C1}
  Let $k$ be a quadratic number field and assume $L/k$ is an
  unramified normal extension with Galois group $\SD_{2^n}=\la
  \sigma,\tau: \sigma^{2m}=\tau^2=1,\tau\sigma\tau=\sigma^{m-1}\ra,$
  with $m=2^{n-2}$. Then $L/\Q$ is normal, and we have
  $\Gal(L/\Q)\simeq \Gamma_1$ with
  $$\Gamma_1 = \la \vr,\sigma,\tau:\vr^2=\tau^2=\sigma^{2m}=1,
     \vr\sigma\vr=\sigma^{-1},\vr\tau\vr=\tau,\tau\sigma\tau=\sigma^{m-1}\ra.$$

   Moreover, there is a unique quartic extension $K/k$ contained in
   $k_\gen$ with\\ $\Cl_2(K)\simeq C_{2^{n-2}}.$
\end{cor}

\begin{proof}
  Since $L=k^2$, we see that $L/\Q$ is normal. The structure of
  $\Gamma$ follows from the previous theorem and
  Proposition~1. Finally, $J=\la \sigma^2\ra$ is the only cyclic
  subgroup of $\SD_{2^n}$ of order $2^{n-2}$. Let $K=L^J$ be the fixed
  field of $J$ in $L$. Then $\Gal(K/k)\simeq \SD_{2^n}/J\simeq
  C_2\times C_2$ and therefore $K\subseteq k_\gen.$
\end{proof}

\begin{cor}\label{C2}
  Let $k$ be a quadratic number field and assume $L/k$ is an
  unramified normal extension with Galois group $\MD_n(2)=\la
  \sigma,\tau: \sigma^{2m}=\tau^2=1,\tau\sigma\tau=\sigma^{m+1}\ra,$
  with $m=2^{n-2}$. Then $L/\Q$ is normal, and we have
  $\Gal(L/\Q)\simeq \Gamma_2$ with
  $$\Gamma_2 = \la \vr,\sigma,\tau:\vr^2=\tau^2=\sigma^{2m}=1,
    \vr\sigma\vr=\sigma^{-1},\vr\tau\vr= \tau,\tau\sigma\tau=\sigma^{m+1}\ra.$$
\end{cor}
The proof follows immediately from the previous theorem and Proposition~1.

\section{Capitulation Gaps} 

For the classification of unramified extensions with semi-dihedral
and modular Galois groups we will use a simple but effective result,
which generalizes an observation made by Couture and Derhem
in \cite{CD} (see also \cite[Prop. 7.17]{LTh}):

\begin{thm}\label{ThCD}
  Let $k$ be a number field whose $2$-class group $\Cl_2(k)$ has
  rank $r \ge 2$. If there is an ideal class in $\Cl(k)[2]$ that
  does not capitulate in any quadratic unramified extension of $k$, then
  there is a unit $\eps \in E_k$ with $\eps \equiv \xi^2 \bmod 4$,
  i.e., such that $k(\sqrt{\eps}\,)/k$ is a quadratic unramified extension.
\end{thm}

The proof will make use of the Selmer group $\Sel_4^+(k)$ (see \cite{LSQ};
this coincides with $R_k$ in \cite{CD} when $k$ is totally complex)
$$ \Sel_4^+(k) = \{\alpha \in k^\times: (\alpha, 2) = (1), \ \alpha \gg 0, \ 
   \alpha \equiv \xi^2 \bmod 4,\ (\alpha) = \fa^2\}/k^{\times\,2} $$
generated by totally positive elements $\alpha \in k^\times$ for which
$k(\sqrt{\alpha}\,)/k$ is unramified at all finite and infinite primes.
The totally positive units congruent to squares modulo $4$ generate a
subgroup $E_4^+/E_k^2$ of $\Sel_4^+(k)$:
$$ E_4^+ = \{\eps \in E_k: \eps \gg 0, \eps \equiv \xi^2 \bmod 4\}. $$
If $\eps \in E_4^+$, then $k(\sqrt{\eps}\,)$ is unramified everywhere.

We will split up the proof of Thm.~\ref{ThCD} into two lemmas both dealing
with the subgroup $\Cl^*(k)$ of ideal classes of order $2$
generated by ideals $\fa$ for which $\fa^2 = (\alpha)$ for some
$\alpha \in \Sel_4^+(k)$. 

\begin{lem}
  We have an exact sequence
   $$ \begin{CD}
    1 @>>> E_4^+/E^2 @>>> \Sel_4^+(k) @>{\phi}>> \Cl^*(k)  @>>> 1.    
  \end{CD} $$
  In particular, $$\#  \Cl^*(k) = \frac{\# \Cl(k)[2]}{(E_4^+:E^2)}. $$
\end{lem}

\begin{proof}
  We consider the map $\phi$ sending an element
  $\alpha k^{\times\,2} \in  \Sel_4^+(k)$ with $(\alpha) = \fa^2$ to the
  ideal class $[\fa] \in \Cl(k)[2]$. Clearly $\phi$ is surjective, and we have
  $\alpha k^{\times\,2} \in \ker \phi$ if and only if
  $(\alpha) = (\beta)^2$, i.e., if and only if $\alpha = \eps \beta^2$.
  Thus $\ker \phi$ is generated by units $\eps \gg 0$ with
  $\eps \equiv \xi^2 \bmod 4$, i.e., $\ker \phi = E_4^+/E^2$.
\end{proof}

\begin{lem}
  Each ideal class in $\Cl^*(k)$ capitulates in some quadratic unramified
  extension of $k$.
\end{lem}

\begin{proof}
  If $[\fa] \in \Cl^*(k)$, then $\fa^2 = (\alpha)$ for some totally positive
  $\alpha \equiv \xi^2 \bmod 4$. The extension $K = k(\sqrt{\alpha}\,)$
  is unramified over $k$, and clearly $\fa \cO_K = (\sqrt{\alpha}\,)$
  capitulates in $K/k$.
\end{proof}

\begin{proof}[Proof of Thm.~\ref{ThCD}]
  An ideal class of order $2$ in $\Cl(k)[2]$ that does not capitulate
  lies in $\Cl(k)[2] \setminus \Cl^*(k)$. This is only possible
  if $(E_4^+:E^2) \ge 2$.
\end{proof}

The condition that there is a class of order $2$ that does not capitulate
in any quadratic unramified extension is group theoretical. Set
$\Gamma = \Gal(k^2/k)$, the Galois group of the second Hilbert $2$-class
field; then $\Gamma/\Gamma' \simeq \Gal(k^1/k) \simeq \Cl_2(k)$, and the
classes of order $2$ correspond to groups $G$ with
$\Gamma' \subset G \subseteq \Gamma$ with $(G:\Gamma') = 2$. The ideal
class $c \simeq G/\Gamma'$ capitulates in some quadratic extension $K/k$
(corresponding to subgroups $H$ of index $2$ in $\Gamma$) if and only if
$G/\Gamma' \subseteq \ker \Ver_{\Gamma \lra H}$,
where $V = \Ver_{\Gamma \lra H}: \Gamma/\Gamma' \lra H/H'$ denotes the
transfer map. We say that a finite $2$-group has a capitulation gap
if there is a subgroup $G$ such that $G/\Gamma'$ is not contained in some
capitulation kernel $ \ker \Ver_{\Gamma \lra H}$.

$$ \begin{array}{ccccccc}
  k      & \hra    & K & \hra    &  k^1    & \hra    & k^2 \\
  \updownarrow & & \updownarrow & & \updownarrow & & \updownarrow \\
  \Gamma & \supset & H & \supset & \Gamma' & \supset & 1
  \end{array} $$

It is clear that abelian $2$-groups do not have a capitulation
gap; in fact if $\Gamma$ is an abelian $2$-group and if $H$ is a
subgroup of index $2$, then $\ker \Ver_{\Gamma \lra H} = \Gamma[2]$
contains all elements of order $2$ and this kernel contains every
subgroup $G$ of $\Gamma$ with order $2$.

\begin{prop}\label{Pro9}
  The semi-dihedral groups $\Gamma = \SD_{2^n}$ and the modular groups
  $\Gamma = \MD_4(n)$ have a capitulation gap.
\end{prop}

\begin{proof}
  In both cases there are three subgroups $H$ of index $2$ in $\Gamma$.
  \begin{enumerate}
  \item $\Gamma \simeq \SD_{2^n}$. The semi-dihedral group of order $4m$
    has the presentation
    $$ \Gamma = \la \sigma, \tau: \sigma^m = -1, \tau^2 = 1,
                    \tau\sigma\tau = \sigma^{m-1} \ra,$$
    and we have $\Gamma' = \la \sigma^2 \ra$.
    The subgroups $H$ of index $2$ in $G$ and their transfer kernels are
    given in the following table:
    $$ \begin{array}{c|ccc}
      \rsp H & \Ver(\sigma)/H' & \Ver(\tau)/H' & \ker \Ver \\ \hline
      \rsp \la \sigma \ra      &         \sigma^{m}H' & H'
           & \la \tau \ra \Gamma' \\ 
      \rsp \la \sigma^2, \tau \ra    &   \sigma^2 H'   & \sigma^{-2}H'
           & \la \sigma\tau \ra \Gamma' \\ 
      \rsp \la \sigma^2, \sigma\tau \ra    & \sigma^2 H'   &  H' 
           & \la \tau \ra  \Gamma'
      \end{array} $$
    In particular, the group $\la \sigma \ra \Gamma'$ is not contained in any
    transfer kernel, hence $\Gamma$ has a capitulation gap.
  \item $\Gamma = \MD_n(2)$: The modular group $\MD_n(2)$ has the presentation
    $$ \Gamma = \la \sigma, \tau: \sigma^{2m} = \tau^2 = 1,
       \tau \sigma \tau = \sigma^{m+1} \ra, $$
    where $m = 2^{n-2}$, and $\Gamma' = \la \sigma^m \ra$.
    The transfer kernels are given in the following table:
     $$ \begin{array}{c|ccc}
      \rsp H & \Ver(\sigma)/H' & \Ver(\tau)/H' & \ker \Ver \\ \hline
      \rsp \la \sigma \ra               & \sigma^{m+2} H' & 1 & 
           \la \tau \ra \Gamma'   \\ 
           \rsp \la \sigma^2, \tau \ra       & \sigma^2 H'     & \sigma^m  &
           \la \tau \sigma^{m/2} \ra \Gamma'
           \\ 
      \rsp \la \sigma^2, \sigma\tau \ra & \sigma^2 H'     & 1  & 
           \la \tau \ra \Gamma' 
      \end{array} $$
  \end{enumerate}
  In particular, the group $\la \sigma^{m/2} \ra \Gamma'$ is not
  contained in any transfer kernel, hence $\Gamma$ has a capitulation gap.
  This completes the proof.
\end{proof}

\begin{cor}
  If $k$ is a complex quadratic number field and if $\Gal(k^2/k)$ has
  a capitulation gap, then $k = \Q(\sqrt{m}\,)$ for some $m \equiv 3 \bmod 4$,
  i.e., the factorization of $\disc k$ into prime discriminants contains
  the factor $-4$.

  In addition, if $\Cl(k)$ has $2$-rank equal to $2$, then $d = -4 pq$
  for positive prime discriminants $p \equiv q \equiv 1 \bmod 4$; in
  this case, the fundamental unit $\eps$ of $\Q(\sqrt{pq}\,)$ has negative
  norm.
\end{cor}

In fact, if $k$ is complex quadratic, the only nontrivial unit is
$\eps = -1$ (the fields with discriminants $-3$ and $-4$ have class
number $1$), and $-1$ is congruent to a square modulo $4$ if and only if
$k(\sqrt{-1}\,)/k$ is an unramified quadratic extension contained in the
genus class field. This happens if and only if $-4$ is a factor in the
factorization of $\disc k$ into prime discriminants.

The second claim follows from the fact that if $N\eps = +1$,
then all classes of order $2$ must capitulate in $K/k$ since
we know that in cyclic unramified extensions $K/k$ we have
$\# \kappa_{K/k} = (K:k)(E_k:N_{K/k}E_K)$. 
Since, in this case, $\eps_{pq}$ remains fundamental in $K$, $-1$
is the norm of a unit only if $N\eps_{pq} = -1$.
  
\subsection*{Calculations with GAP}

We now use GAP to find all $2$-groups of small order
with a capitulation gap.

\begin{prop}
  The only $2$-groups of order $16$ with a capitulation gap are the
  following:
   \begin{center} \begin{tabular}{l|cc}
       \rsp GAP & $16.06$  & $16.08$ \\   \hline 
       \rsp HS  & $16.011$ & $16.013$ \\
       \rsp     & $\MD_4(2)$ & $\SD_{16}$
   \end{tabular} \end{center}
\end{prop}

For groups of order $32$ we find:

\begin{prop}
  The only $2$-groups of order $32$ with a capitulation gap are
  the following:

  $$ \begin{array}{l|ccccc}
  \rsp GAP & 32.04 & 32.13 & 32.15 & 32.17 & 32.19 \\
      \hline 
  \rsp HS  & 32.019 & 32.030 & 32.032 & 32.022 & 32.050 \\
  \rsp \Gamma' & [2] & [4] & [4] & [2] & [8] \\
  \rsp \Gamma/\Gamma' & [4, 4] & [2, 4]  & [2, 4]  & [2, 8]  & [2, 2] 
  \end{array} $$
  All these groups have trivial Schur multiplier. 
  The groups with GAP ID 32.04 and 32.13 have two capitulation gaps.
\end{prop}

In fact, the group $\Gamma = 32.04$ satisfies
$\Gamma/\Gamma' \simeq C_4 \times C_4$, and in each of the three
quadratic extensions, the same ideal class capitulates. This is
also true for the group  $\Gamma = 32.13$ with
$\Gamma/\Gamma' \simeq C_2 \times C_4$.

Here the group with GAP ID 32.17 is the modular group $\MD_8(2)$ and
32.19 is the semidihedral group $\SD_{32}$.

\begin{prop}
  The $2$-groups of order $64$ with two capitulation gaps are
  the following:
    $$ \begin{array}{c|ccccccc}
    \rsp \Gamma & 64.03 & 64.14 & 64.15 & 64.27
                & 64.28 & 64.46 & 64.48  \\ \hline 
    \rsp \Gamma'        & [2]   & [2,4] & [4]   & [2]   & [4]   & [8] & [8] \\
    \rsp \Gamma/\Gamma' & [4,8] & [2,4] & [2,8] & [4,8]
                                            & [4,4] & [2,4] & [2,4]       
  \end{array} $$

  Those with a single capitulation gap are
  
  $$ \begin{array}{c|cccccc}
   \rsp \Gamma
      & 64.09 & 64.10 & 64.45 & 64.49 & 64.51 & 64.53 \\ \hline 
    \rsp \Gamma'        & [2,4] & [2,4] & [4]   & [8]   & [2]    & [16] \\   
    \rsp \Gamma/\Gamma' & [2,4] & [2,4] & [2,8] & [2,4] & [2,16] & [2,2]    
  \end{array} $$
  These groups have trivial Schur multiplier with the exception of
  $\cM(64.09) = [2, 2, 2]$ and $\cM(64.10) =  \cM(64.14) =  [4]$.
  The groups with GAP IDs $64.51$ and $64.53$ are the modular and the
  semidihedral group of order $64$, respectively.
\end{prop}

\section{A few results from class field theory}

The following proposition collects several classical results by
R\'edei, Reichardt, and Scholz on quartic cyclic unramified extensions:

\begin{prop}\label{Pr2}
  The real quadratic number field $k$ with discriminant $d$ admits a cyclic
  quartic extension $L/k$ unramified at all finite primes if and only if
  $d$ has a $C_4$-factorization, i.e., if and only if
  $d = d_1d_2$ for discriminants $d_1, d_2$ such that
  $(d_1/p_2) = (d_2/p_1) = +1$ for all $p_1 \mid d_1$ and all $p_2 \mid d_2$.

  This extension is unramified at infinity if and only if $(\eps_1/d_2) = +1$,
  where $\eps_1$ is the fundamental unit of $k_1 = \Q(\sqrt{d_1}\,)$, which by
  Scholz's reciprocity law is equivalent to
  $$ \Big(\frac{\eps_1}{d_2}\Big) =
     \Big(\frac{d_1}{d_2}\Big)_4 \Big(\frac{d_2}{d_1}\Big)_4. $$ 
\end{prop}

Now we claim:

\begin{lem}\label{L2}
  Assume that $d_1, d_2$ are positive prime discriminants with
  $(d_1/d_2) = +1$. If $F = \Q(\sqrt{d_1d_2}\,)$ does not have an
  everywhere unramified $C_4$-extension, then $(d_1/d_2)_4 (d_2/d_1)_4 = -1$.
\end{lem}

This is clear since if $(d_1/d_2)_4 (d_2/d_1)_4 = 1$, then the cyclic
quartic extension of $F$ corresponding to the $C_4$-factorization
$d = d_1 \cdot d_2$ is also unramified at the infinite primes.

\section{Unramified $\SD_{16}$-extensions of complex quadratic number fields}

It is known (see \cite{LTh} and \cite{BLSI}) for which complex
quadratic number fields $k$ the Galois group of the $2$-class field
tower $k^2/k$ is semi-dihedral:

\begin{thm}
  Let $k$ be a complex quadratic number field with discriminant $d$.
  Then the following assertions are equivalent.
  \begin{enumerate}
  \item[(i)] The $2$-class field tower of $k$ has Galois group $\SD_{2^n}$,
  \item[(ii)] We have $\disc k = -4pq$ for primes $p$ and $q$ with 
    $p \equiv 1 \bmod 8$, $q \equiv 5 \bmod 8$, and $(p/q) = -1$.
  \end{enumerate}
  In this case, $L = k^2$ is cyclic over $k_q = k(\sqrt{q}\,)$, with
  $\Cl_2(k_q) = [2^{n-1}]$.
\end{thm}

We will now give a new proof of this result using Theorem~\ref{ThCD}.

\begin{proof}
  Assume first that the $2$-class field tower of the complex quadratic
  number field $k$ has Galois group $G = \SD_{2^n}$. Since $G/G' \simeq (2,2)$,
  the discriminant $d = \disc k$ of $k$ is a product of three prime
  discriminants.   Prop.~\ref{Pro9} and its corollary imply that 
  $d = -4d_2d_3$ for prime discriminants $d_2$ and $d_3$. If $d_2 < 0$,
  then $\Q(i,\sqrt{d_2d_3}\,)$
  has class number $2$, and the class field tower is abelian. Thus
  $d_2 = p$ and $d_3 = q$ for primes $p \equiv q \equiv 1 \bmod 4$.

  If $p \equiv q \bmod 8$, then $d = -4 \cdot pq$ is a $C_4$-factorization,
  contradicting the fact that $\Cl_2(k) \simeq [2,2]$. Thus, say,
  $p \equiv 1 \bmod 8$ and $q \equiv 5 \bmod 8$.

  If $(p/q) = +1$, then $d = -4q \cdot p$ is a $C_4$-factorization:
  contradiction. Thus $(p/q) = -1$.

  We have shown that if the $2$-class field tower of the complex
  quadratic number field $k$ has Galois group $G = \SD_{2^n}$, then $d
  = -4pq$ for primes $p \equiv 1 \bmod 8$ and $q \equiv 5 \bmod 8$
  with $(p/q) = -1$.
  
  Conversely, assume that $d = -4pq$ for primes $p \equiv 1 \bmod 8$ and
  $q \equiv 5 \bmod 8$ with $(p/q) = -1$. Then $\Cl_2(k) \simeq (2,2)$,
  hence the Galois group of the $2$-class field tower of $k$ is either
  $(2,2)$, dihedral, quaternion or semi-dihedral. The only group with
  a capitulation gap among them is the semi-dihedral group, so we only
  need to show that there is an ideal class of order $2$ in $k$ that does
  not capitulate in any of the three quadratic unramified extensions.

  We begin by observing that $\Cl_2(k) \simeq [2,2]$, and the ideals
  whose classes generate the $2$-class groups are the prime ideal
  $\fp_2 = (2,1+\sqrt{-pq}\,)$ above $2$ and the prime ideal
  $\fp_p = (p,\sqrt{-pq}\,)$ above $p$.

  Recall that the order of the subgroup $\kappa_{K/k}$ of $\Cl_2(k)$
  consisting of ideal classes that capitulate in a quadratic
  unramified extension $K/k$ is given by $|\kappa_{K/k}| = 2(E_k:N_{K/k}E_K).$

  Since $N(\eps_{pq}) = -1$, only one nontrivial ideal class capitulates
  in $K = k(i)$, and this is the ideal class generated by the prime above $2$.
  Similarly, $|\kappa| = 2$ for the other two quadratic unramified
  extensions since $\eps_p$ and $\eps_q$ have norm $-1$. But since
  $\fp_p\cO_K = (\sqrt{p}\,)$ in $K = k(\sqrt{p}\,)$ and 
  $\fp_p\cO_K = (\sqrt{-p}\,)$ in $K = k(\sqrt{q}\,)$, the ideal class
  generated by $\fp_2\fp_p$ does not capitulate in any quadratic
  unramified extension of $k$. This completes the proof.
\end{proof}

The Galois group of the $2$-class field tower of the quadratic number
field $\Q(\sqrt{-pq}\,)$ with $p \equiv 1 \bmod 8$, $q \equiv 5 \bmod 8$
and $(p/q) = -1$ is $\SD_{16}$ if and only if the class number $h$ of
$\Q(\sqrt{-p}\,)$ is $h \equiv 4 \bmod 8$. It is well known that
$h \equiv 0 \bmod 8$ if and only if $(\frac{1+i}{a+bi}) = 1$, where
$p = a^2 + b^2$, and that this is equivalent to $p = x^2 + 32y^2$
(see, for example, \cite{BC}). Thus $h \equiv 4 \bmod 8$ if and only
if $p = 4x^2 + 4xy + 9y^2$.

\begin{figure}[ht!]  
\begin{tikzpicture}[node distance=1.5cm]
  \node (G) at (0,0) {$\la \sigma, \tau \ra$};
  \node (C8) [above of=G] {$\la \sigma \ra$};
  \node (H8) [right of=C8] {$\la \sigma^2, \sigma\tau \ra$};
  \node (D8) [left  of=C8] {$\la \sigma^2, \tau \ra$};
  \draw (G) -- (C8);
  \draw (G) -- (H8);
  \draw (G) -- (D8);
  \node (C4) [above of=C8] {$\la \sigma^2 \ra$};
  \node (C2) [above of=C4] {$\la \sigma^4 \ra$};
  \node (E)  [above of=C2] {$\la 1 \ra$};
  \draw (C8) -- (C4);
  \draw (C4) -- (C2);
  \draw (C2) -- (E);
  \draw (D8) -- (C4);
  \draw (H8) -- (C4);
  \node (k1)  [above of=D8] {$\la \sigma^4, \tau \sigma^2\ra$};
  \node (O)   [left  of=k1] {};
  \node (k1') [left  of=O] {$\la \sigma^4, \tau \ra$};  
  \node (k2)  [above of=H8] {$\la \tau\sigma \ra$};
  \node (k2') [right of=k2] {$\la \tau\sigma^3 \ra$};
  \draw (D8) -- (k1);
  \draw (D8) -- (k1');
  \draw (H8) -- (k2);
  \draw (H8) -- (k2');
  \draw (k1) -- (C2);
  \draw (k1') -- (C2);
  \draw (k2) -- (C2);
  \draw (k2') -- (C2);
  \node (K1)  [above of=k1] {$\la \tau\sigma^4 \ra$};
  \node (K1') [left  of=K1] {$\la \tau \ra$};
  \node (K2)  [above of=k1'] {$\la \tau\sigma^2 \ra$};
  \node (K2') [left  of=K2] {$\la \tau\sigma^6 \ra$};
  \draw [line width=6pt,draw=white] (k1) -- (K1); 
  \draw [line width=6pt,draw=white] (k1) -- (K1');
  \draw (k1) -- (K1);
  \draw (k1) -- (K1');
  \draw (k1') -- (K2);  
  \draw (k1') -- (K2');
  \draw (K1) -- (E);
  \draw (K1') -- (E);
  \draw (K2) -- (E);
  \draw (K2') -- (E);
\end{tikzpicture}
\caption{Subgroups of $\SD_{16}$}
\end{figure}
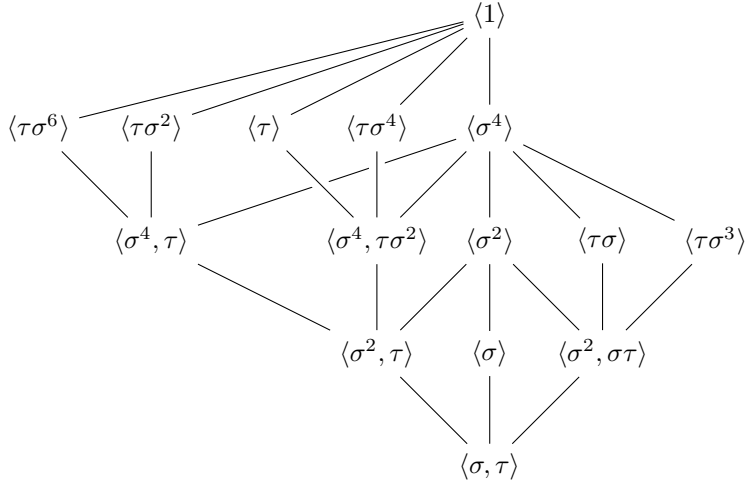

$\SD_{16}$ has three subgroups $G_i$ of index $2$:
$$ \la \sigma \ra \simeq C_8, \quad 
   \la \sigma^2, \sigma\tau \ra \simeq H_8, \quad \text{and} \quad 
   \la \sigma^2, \tau \ra \simeq D_4. $$
These groups of index $2$ fix the three quadratic subextensions of $L/k$;
here $L$ is cyclic over $k_q = k(\sqrt{q}\,)$, a dihedral extension of
$k_i = k(i)$, and a quaternion extension of $k_p = k(\sqrt{p}\,)$.

\subsection*{Construction of unramified $\SD_{16}$-extensions}

The construction of the lower pieces of the $2$-class field tower
of these fields is well known. The genus class field of $k = \Q(\sqrt{-pq}\,)$
is $k_\gen = \Q(i,\sqrt{p},\sqrt{q}\,)$, and if we set $p = a^2 + b^2$
with $a \equiv 1 \bmod 4$, then $E = k_\gen(\sqrt{a+bi}\,)$ is an unramified
quadratic extension of $K = k_\gen$ with $\Gal(E/k) \simeq D_4$.

The two conjugate quadratic unramified extensions of $k_p = k(\sqrt{p}\,)$
are given by $K_p = k(\sqrt{\eps_p}\,)$ and $K_p' = k(\sqrt{\eps_p'}\,)$,
where $\eps_p > 1$ is the fundamental unit of $\Q(\sqrt{p}\,)$ and $\eps_p'$
its conjugate. The well known fact that $(a+bi)\eps_p$ is a square in $k_\gen$
is the source of some reciprocity laws related to Scholz's reciprocity law.

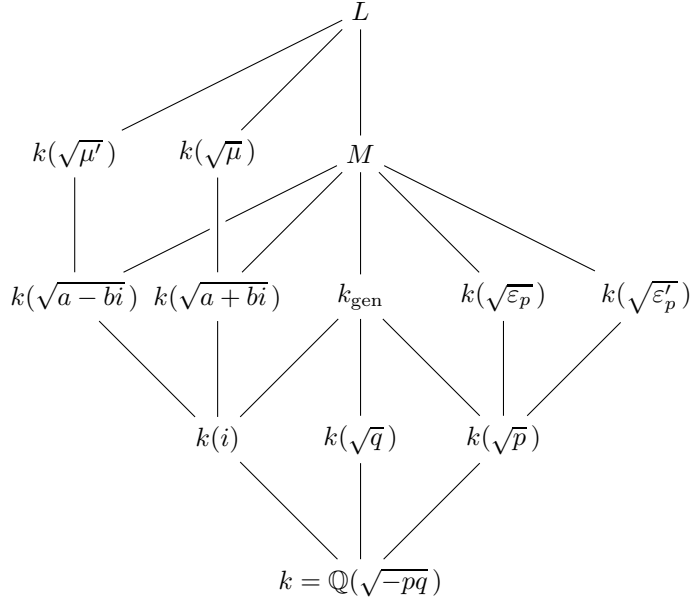
\begin{figure}[ht!]  
\begin{tikzpicture}[node distance=1.9cm]
  \node (G) at (0,0) {$k = \Q(\sqrt{-pq}\,)$};
  \node (C8) [above of=G] {$k(\sqrt{q}\,)$};
  \node (H8) [right of=C8] {$k(\sqrt{p}\,)$};
  \node (D8) [left  of=C8] {$k(i)$};
  \draw (G) -- (C8);
  \draw (G) -- (H8);
  \draw (G) -- (D8);
  \node (C4) [above of=C8] {$k_\gen$};
  \node (C2) [above of=C4] {$M$};
  \node (E)  [above of=C2] {$L$};
  \draw (C8) -- (C4);
  \draw (C4) -- (C2);
  \draw (C2) -- (E);
  \draw (D8) -- (C4);
  \draw (H8) -- (C4);
  \node (k1)  [above of=D8] {$k(\sqrt{a+bi}\,)$};
  \node (k1') [left  of=k1] {$k(\sqrt{a-bi}\,)$};
  \node (k2)  [above of=H8] {$k(\sqrt{\eps_p}\,)$};
  \node (k2') [right of=k2] {$k(\sqrt{\eps_p'}\,)$};
  \draw (D8) -- (k1);
  \draw (D8) -- (k1');
  \draw (H8) -- (k2);
  \draw (H8) -- (k2');
  \draw (k1) -- (C2);
  \draw (k1') -- (C2);
  \draw (k2) -- (C2);
  \draw (k2') -- (C2);
  \node (K1)  [above of=k1]  {$k(\sqrt{\mu}\,)$};
  \node (K2)  [above of=k1'] {$k(\sqrt{\mu'}\,)$};
  \draw [line width=6pt,draw=white] (k1) -- (K1); 
  \draw (k1) -- (K1); 
  \draw (k1') -- (K2);
  \draw (K1) -- (E);
  \draw (K2) -- (E);
\end{tikzpicture}
\caption{Subextensions of the $2$-class field tower of $\Q(\sqrt{-pq}\,)$
  (with two conjugate fields omitted)}\label{ASD16}
\end{figure}

We will need the following result from \cite[Thm. 8]{LemD}:

\begin{prop}\label{PropT}
  Let $\lambda = a+bi \equiv 1 \bmod 4$ be a Gaussian integer with
  norm $m$. If the equation
  $$ A^2 - \lambda B^2 - \ov{\lambda} C^2 = 0 $$
  has nontrivial solutions in Gaussian integers, then the solution
  $(\alpha, \beta, \gamma)$ can be chosen in such a way that
  $\mu = \alpha + \beta \sqrt{\lambda} \equiv 1 \bmod 2+2i$. In this
  case, $L = \Q(i,\sqrt{\mu}\,)$ is a cyclic extension of degree $8$ over
  $\Q(\sqrt{-m}\,)$ unramified outside $2$. The extension $L/k$ is
  unramified everywhere if and only if $(\frac2{a+b}) = +1$.
\end{prop}

Now we claim:

\begin{prop}
  Let $p \equiv 1 \bmod 8$ and $q \equiv 5 \bmod 8$ be primes
  with $(\frac pq) = -1$, and write $\pi = a+bi$, where
  $p = a^2 + b^2$ and $\pi \equiv 1 \bmod 4$. Then the equation
  \begin{equation}\label{ESD1} A^2 - \pi B^2 = \pi' \Gamma^2  \end{equation}
  has a nontrivial solution $(\alpha, \beta, \gamma)$ in $\Z[i]$.

  Setting $\mu = \alpha + \beta \sqrt{\pi}$, the extension
  $K = k(\sqrt{p}, \sqrt{\pi}, \sqrt{\mu}\,)$ is a $C_8$-extension
  of $k$. We can choose $\alpha$ and $\beta$ in such a way that
  $\mu \equiv 1 \bmod 2+2i$, and then $K/k$ is
  \begin{itemize}
  \item unramified if $p = x^2 + 32y^2$,
  \item ramified at $2$ if $p = 4x^2 + 4xy + 9y^2$. In this case,
    $L =  k(\sqrt{p}, \sqrt{\pi}, \sqrt{\rho\mu}\,)$ is an unramified
    $\SD_{16}$-extension of $k$, where $\rho = c + di \equiv 1 \bmod 2+2i$
    and $q = c^2 + d^2$.
  \end{itemize}
\end{prop}

Assume now that $p \equiv 1 \bmod 8$, $q \equiv 5 \bmod 8$
and $(p/q) = -1$, and write $p = a^2 + b^2$ with $b \equiv 0 \bmod 4$.
In \cite{LemD} we have shown that the equation
$$ \alpha^2 - (a+bi)\beta^2 = (a-bi)\gamma^2 $$
is nontrivially solvable with $\alpha, \beta, \gamma \in \Z[i]$, and
that the solution can be chosen in such a way that
$\mu = \alpha + \beta \sqrt{a+bi} \equiv 1 \bmod 2+2i$. Thus 
$\mu \equiv 1$ or $\equiv 3+2i \bmod 4$, and $k_\gen(\sqrt{\mu}\,)$ is an
unramified cyclic extension of degree $4$ if and only if $p = x^2 + 32y^2$.

If  $p \ne x^2 + 32y^2$ then clearly $\mu(c+di)$ is congruent to a square
modulo $4$, where $q = c^2 + d^2$ and $d$ is even. This implies that
$L = M(\sqrt{\mu(c+di)}\,)$ is a quadratic unramified extension, which
therefore must coincide with the second Hilbert $2$-class field $k^2$ of $k$.
The fact that $\Gal(L/k) \simeq \SD_{16}$ can also be checked directly using
the methods presented in \cite{LemD}.

\begin{figure}[ht!]  
\begin{tikzpicture}[node distance=1.5cm]
  \node (G) at (0,0) {$[2,2]$};
  \node (C8) [above of=G] {$[8]$};
  \node (H8) [right of=C8] {$[2,2]$};
  \node (D8) [left  of=C8] {$[2,2]$};
  \draw (G) -- (C8);
  \draw (G) -- (H8);
  \draw (G) -- (D8);
  \node (C4) [above of=C8] {$[4]$};
  \node (C2) [above of=C4] {$[2]$};
  \node (E)  [above of=C2] {$[1]$};
  \draw (C8) -- (C4);
  \draw (C4) -- (C2);
  \draw (C2) -- (E);
  \draw (D8) -- (C4);
  \draw (H8) -- (C4);
  \node (k1)  [above of=D8] {$[2,2]$};
  \node (k1') [left  of=k1] {$[2,2]$};
  \node (k2)  [above of=H8] {$[4]$};
  \node (k2') [right of=k2] {$[4]$};
  \draw (D8) -- (k1);
  \draw (D8) -- (k1');
  \draw (H8) -- (k2);
  \draw (H8) -- (k2');
  \draw (k1) -- (C2);
  \draw (k1') -- (C2);
  \draw (k2) -- (C2);
  \draw (k2') -- (C2);
  \node (K1)  [above of=k1]  {$[2]$};
  \node (K2)  [above of=k1'] {$[2]$};
  \draw [line width=6pt,draw=white] (k1) -- (K1); 
  \draw (k1) -- (K1); 
  \draw (k1') -- (K2);
  \draw (K1) -- (E);
  \draw (K2) -- (E);
\end{tikzpicture}
\caption{$2$-class groups of the subextensions of $L/k$ for unramified
  $\SD_{16}$-extensions of $k = \Q(\sqrt{-pq}\,)$}
\end{figure}

\begin{table}
  $$ \begin{array}{r|r|rr|rr|c}
  \rsp   p &  q &   a &  b & c & d & \mu\rho  \\ \hline
  \rsp  17 &  5 &   1 &  4 & 1 & 2 & (2+2i + \sqrt{1+4i}\,)(1+2i) \\
  \rsp  73 &  5 &  -3 &  8 & 1 & 2 & (2+2i + (3+2i)\sqrt{-3+8i}\,)(1+2i) \\
  \rsp  97 &  5 &   9 &  4 & 1 & 2 & (2+2i + \sqrt{9+4i}\,)(1+2i) \\
  \rsp  73 & 13 &  -3 &  8 & 3 & 2 & (2+2i + (3+2i)\sqrt{-3+8i}\,)(3+2i)  \\
  \rsp 193 &  5 &  -7 & 12 & 1 & 2 & (2+2i + (7+4i)\sqrt{-7+12i}\,)(1+2i) \\
  \rsp  89 & 13 &   5 &  8 & 3 & 2 & (10+10i) + (3+2i)\sqrt{5+8i}\,)(3+2i) \\
  \rsp 233 &  5 &  13 &  8 & 1 & 2 & (14+14i) + (3+2i)\sqrt{13+8i}\,)(1+2i)\\
  \rsp  97 & 13 &   9 &  4 & 3 & 2 & (2+2i + \sqrt{9+4i}\,)(3+2i) \\
  \rsp  73 & 29 &  -3 &  8 & 5 & 2 & (2+2i + (3+2i)\sqrt{-3+8i}\,)(5+2i)  \\
  \rsp 433 &  5 &  17 & 12 & 1 & 2 & (12 + 12i + (3+2i)\sqrt{17-12i}\,)(1+2i) \\
  \rsp 193 & 13 &  -7 & 12 & 3 & 2 & (2+2i + (7+4i)\sqrt{-7+12i}\,)(3+2i) \\
  \rsp  89 & 29 &   5 &  8 & 5 & 2 & (10+10i) + (3+2i)\sqrt{5+8i}\,)(5+2i) \\
  \rsp  97 & 29 &   9 &  4 & 5 & 2 & (2+2i + \sqrt{9+4i}\,)(5+2i) 
  \end{array} $$
  \caption{Unramified $\SD_{16}$-extensions of $k = \Q(\sqrt{-pq}\,)$;
  here  $K = \Q(i,\sqrt{p},\sqrt{q}\,)$, $E = K(\sqrt{a+bi}\,)$
  and $L = E(\sqrt{\mu\rho}\,)$.}\label{TabSD}
\end{table}


The signs of $b$ and $d$ can be chosen arbitrarily; the choice of signs
in $\mu\rho$ is determined by the fact that it must give rise to a
solution of (\ref{ESD1}). For example we have
$$ (2+2i)^2 - (-3+8i)(3+2i)^2 = (2+3i)^2(-3-8i). $$

\section{Unramified $\MD_4(2)$-extensions of complex quadratic number fields}

In the second part of his report on class field theory \cite[p. 173]{Hasse},
Helmut Hasse asked whether in a subextension $K/k$ of the Hilbert
class field $k^1/k$ of a number field there is always a subgroup
of $\kappa_{K/k}$, the group of ideals in $k$ becoming principal in $K$,
that is isomorphic to $\Gal(K/k)$. Then he presents
a group theoretical counterexample due to Furtw\"angler:
If there is a number field $k$ such that $\Gal(k^2/k) \simeq \MD_4(2)$,
then in the two cyclic quartic subextensions of $k^1/k$, the group of ideals
that capitulate is $\simeq (2,2)$, and the ideals of order $4$ in the
class group only capitulate in $k^1$. Artin then showed that this
situation is realized by $k = \Q(\sqrt{-65}\,)$.

\begin{figure}[ht!]
\begin{tikzpicture}[node distance=1.5cm]
  \node (G) at (0,0) {$\la \sigma, \tau \ra$};
  \node (C8)  [above of =G] {$\la \sigma\tau \ra$};
  \node (C81) [left  of=C8] {$\la \sigma \ra$};
  \node (C24) [right of=C8] {$\la \sigma^2, \tau \ra$};
  \draw (G) -- (C8);
  \draw (G) -- (C81);
  \draw (G) -- (C24);
  \node (C4) [above of=C8] {$\la \sigma^2 \ra $};
  \draw (C8) -- (C4);
  \draw (C81) -- (C4);
  \draw (C24) -- (C4);
  \node (H) [above of=C24] {$\la  \sigma^2\tau \ra $};
  \node (H2) [right of=H]  {$\la  \sigma^4, \tau \ra $};
  \draw (C24) -- (H);
  \draw (C24) -- (H2);
  \node (C2) [above of=H] {$\la  \sigma^4 \ra $};
  \node (C21) [above of=H2] {$\la \sigma^4\tau \ra $};
  \node (C22) [right of=C21] {$\la \tau \ra $};
  \draw (C4) -- (C2);
  \draw (H) -- (C2);
  \draw (H2) -- (C2);
  \draw (H2) -- (C21);
  \draw (H2) -- (C22);
  \node (E) [above of=C21] {$\la 1 \ra $};
  \draw (C2) -- (E);
  \draw (C21) -- (E);
  \draw (C22) -- (E);
\end{tikzpicture}
\caption{Subgroups of $\MD_4(2)$}
\end{figure}
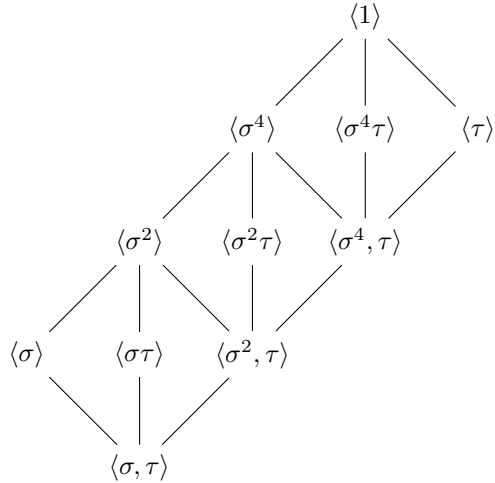

In this section we will explicitly construct unramified extensions
$L/k$ of complex quadratic number fields $k$ with Galois group
$\Gal(L/k) \simeq \MD_4(2)$, the modular group of order $16$.
More generally, the group $\MD_n(2)$ is a semidirect product
$\MD_n(2) \simeq C_{2^{n-1}} \rtimes C_2$, where $C_2 = \la \tau \ra$
acts on $ C_{2^{n-1}} = \la \sigma \ra$ as $\tau \sigma \tau = \sigma^{2^{n-2}+1}$.
In particular,
$$\MD_n(2) = \la \sigma, \tau: \sigma^{2^{n-1}} = \tau^2 = 1,
              \tau \sigma \tau = \sigma^{2^{n-2}+1}. \ra $$

The complex quadratic number fields whose $2$-class field tower              
has Galois group $\MD_n(2)$ were  classified in \cite{BLSI}:

\begin{thm}\label{TClassM}
  Let $k$ be a complex quadratic number field with discriminant $d$.
  Then the following assertions are equivalent.
  \begin{enumerate}
  \item[(i)] The $2$-class field tower of $k$ has Galois group $\MD_n(2)$
    for some $n = 2^m$ with $m \ge 4$;
  \item[(ii)] $\disc k = -4pq$ for primes $p$ and $q$ with 
    $p \equiv q \equiv 5 \bmod 8$ and $(p/q) = -1$.
  \end{enumerate}
  The integer $n$ is determined by the $2$-class number of $k$; in fact
  we have $\Cl_2(k) \simeq [2^{n-2},2]$.
\end{thm}

Thus $\Gal(k^2/k) \simeq \MD_4(2)$ if and only if $\Cl_2(k) \simeq [4,2]$.
Kaplan \cite[Thm. $B_1'$]{Kaplan} has proved the following

\begin{prop}
  Let $k = \Q(\sqrt{-pq}\,)$ for primes $p \equiv q \equiv 5 \bmod 8$
  with $(\frac pq) = -1$. Then $\Cl_2(k) \simeq [2,4]$ if and only if
  $$ \Big(\frac{pq}2\Big)_4  \Big(\frac{2p}{q}\Big)_4 \Big(\frac{2q}{p}\Big)_4
     = -1, $$
  and  $\Cl_2(k) \simeq [2,2^n]$ for some $n \ge 3$ otherwise.
\end{prop}

Observe that $(\frac{a}2)_4 = (-1)^{\frac{a-1}8}$ for integers
$a \equiv 1 \bmod 8$ by definition.

Let us also sketch a proof of Theorem~\ref{TClassM} using Theorem~\ref{ThCD}.
We begin by assuming that $\Gal(k^2/k) \simeq M_n(2)$, where $n = 2^m$.
As before, we must have $d = -4pq$ for primes $p \equiv q \equiv 1 \bmod 4$.
If $(p/q) = +1$, then $\Q(\sqrt{pq}\,)$ admits a $C_4$-extension
unramified at all finite primes, hence $k$ admits an unramified
$D_4$-extension, which is impossible. Thus $(p/q) = -1$. Since
$\Cl_2(k) \simeq [2^{n-2},2]$, there must be a unique $C_4$-factorization
of $d$. Since $(\frac pq) = -1$ we conclude that the
$C_4$-factorization is $d = -4 \cdot pq$; then $(\frac{pq}2) = +1$
implies $p \equiv q \bmod 8$.

If $p \equiv 1 \bmod 8$, then $\Q(\sqrt{-p}\,)$ admits an unramified
$C_4$-extension, which would imply as above that $k$ has an unramified
$D_4$-extension; thus $p \equiv q \equiv 5 \bmod 8$.

Now assume that $d = -4pq$ for primes $p \equiv q \equiv 5 \bmod 8$
with $(\frac pq) = -1$. We know that $\Cl_2(k) \simeq [2^m,2]$ for some
integer $m \ge 2$. We will sketch the proof that $\Gal(k^2/k) \simeq M_n(2)$:
By the ambiguous class number formula, $\Cl_2(k(\sqrt{q}\,))$ is cyclic,
and Kuroda's class number formula implies that it has order $h_2(k)$.
Thus $h_2(k^1) = 2$.

When $\Cl_2(k) = [4^*,2]$ and $h_2(k^1) = 2$, then $\Gal(k^2/k)$ is either
metacyclic, modular, or non-metacyclic. For the nonmodular cases, there's
a capitulation kernel of order $4$. Hence modular is the only one that
could have gaps.

\medskip

Now write $p = a^2 + b^2$ and $q = c^2 + d^2$ with $a$ and $c$ odd, and set
$\pi = a+bi$ and $\rho = c+di$. Then $k(i,\sqrt{\pi\rho}\,)$ and
$k(i,\sqrt{\pi\rho'}\,)$ are the quadratic unramified extensions of
$k_i = k(i)$ different from $k_\gen = \Q(i,\sqrt{p},\sqrt{q}\,)$.
In the following we will need the quadratic residue symbol
$[\frac{\cdot}{\cdot}]$ in the Gaussian integers:

\begin{prop}
  Let $p \equiv q \equiv 5 \bmod 8$ be primes with $(p/q) = -1$,
  and assume that $k = \Q(\sqrt{-pq}\,)$ has $2$-class group $[4,2]$.
  Fix a prime $\pi = a+bi$ dividing $\pi$ and choose
  $\rho = c+di$ dividing $q$ in such a way that $[\frac{\pi}{\rho}] = +1$.
  Then $F = k_i(\sqrt{\pi\rho}\,)$ has $2$-class group $\Cl_2(F) \simeq [2,2]$,
  and $F' = k_i(\sqrt{\pi\rho'}\,)$ has $2$-class group $\Cl_2(F') \simeq [4]$.
\end{prop}

\begin{figure}[ht!]  
\begin{tikzpicture}[node distance=1.5cm]
  \node (Qi) at (0,0) {$\Q(i)$};
  \node (K) [above of=Qi] {};
  \node (Qipq) [left  of=K] {$\Q(i,\sqrt{pq}\,)$};
  \node (E) [right  of=K] {$E' = \Q(i,\sqrt{\pi\rho'}\,)$};
  \node (F) [above of=K] {$F'=\Q(i,\sqrt{pq},\sqrt{\pi\rho'}\,)$};
  \draw (Qi) -- (Qipq);
  \draw (Qi) -- (E);
  \draw (Qipq) -- (F);
  \draw (E) -- (F);
\end{tikzpicture}
\caption{The extension $F'/\Q(i,\sqrt{pq}\,)$.}
\end{figure}
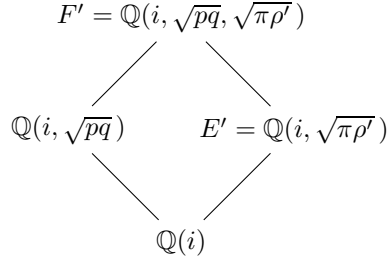

\begin{proof}
  Since $\Gal(L/k_\gen) \simeq C_4$, one of $\Gal(L/F)$ and $\Gal(L/F')$
  is cyclic, and the other one is bicyclic. Thus it is sufficent to
  prove that $\Cl_2(F')$ is cyclic. This is done by applying the
  ambiguous class number formula to $F'/E'$, where
  $E' =  \Q(i, \sqrt{\pi\rho}\,)$.

  Observe that
  $[\frac{\pi\rho'}{\rho}] = [\frac{\pi}{\rho}][\frac{\rho'}{\rho}]
  = [\frac{\rho'}{\rho}]$; but $\rho' = c-di \equiv 2c \bmod \rho$
  shows that
  $[\frac{\rho'}{\rho}] = (\frac{2c}{q}) = (\frac{2}{q}) (\frac{c}{q}) = -1$
  since $q \equiv 5 \bmod 8$ and
  $(\frac{c}{q}) = (\frac{q}{c}) = (\frac{c^2+d^2}{c}) = (\frac{d}{c})^2 = +1$.
  This implies that $\rho$ is inert in $E'/\Q(i)$. Similarly, 
  $[\frac{\pi'\rho}{\pi}] = -1$ implies that $\pi$ is inert in $E'/\Q(i)$.

  We know that $E'$ has odd class number; since $\pi'$ and $\rho'$ are inert
  in $E'/\Q(i)$, exactly two primes ramify in $F'/E'$, and the ambiguous
  class number formula implies that $\# \Am_2(F'/E') \le 2$.
  Thus $\Cl_2(F')$ is cyclic. 
\end{proof}

\subsection*{Construction of unramified $\MD_4(2)$-extensions}
We now turn to the explicit construction of unramified $\MD_4(2)$-extensions
of $k$.

\begin{prop}
  Let $p \equiv q \equiv 5 \bmod 8$ be as above.
  If $[\pi/\rho] = +1$, then the equation
  \begin{equation}\label{Epirho}
    A^2 - \pi \rho B^2 - \pi'\rho' \Gamma^2 = 0
  \end{equation}  
  has a nontrivial solution $(\alpha, \beta, \gamma)$ in $\Z[i]$.

  Setting $\mu = \alpha + \beta \sqrt{\pi\rho}$, the extension
  $K = k(\sqrt{pq}, \sqrt{\pi \rho}, \sqrt{\mu}\,)$ is a $C_8$-extension
  of $k$, which is ramified at $2$.

  In addition, $\alpha$ and $\beta$ can be chosen in such a way that
  $L = k(\sqrt{pq}, \sqrt{\pi \rho}, \sqrt{\pi \mu}\,)$ is an unramified
  $\MD_4(2)$-extension of $k$.
\end{prop}

\begin{proof}
  For proving the solvability of the equation we have to verify the
  solvability everywhere locally. We already know that 
  $[\frac{\pi \rho}{\pi'}] = +1$ and $[\frac{\pi \rho}{\rho'}] = +1$.
  Thus it remains to verify solvability at the prime $1+i$ above $2$. But this
  taken care of by the product formula for the Hilbert symbol.

  By Prop.~\ref{PropT}, $\mu \equiv 1 \bmod 2+2i$ is not a square modulo $4$,
  hence $\pi \mu \equiv \pm 1 = i^2 \bmod 4$. This implies that $L/k$ is
  unramified everywhere.
\end{proof}

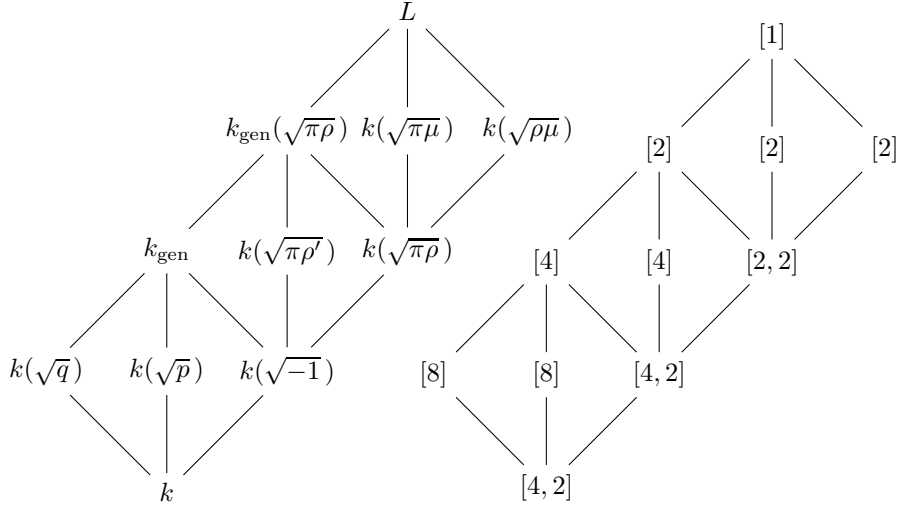
\begin{figure}[ht!]             
\begin{tikzpicture}[node distance=1.6cm]
  \node (G) at (0,0) {$k$};
  \node (C8)  [above of =G] {$k(\sqrt{p}\,)$};
  \node (C81) [left  of=C8] {$k(\sqrt{q}\,)$};
  \node (C24) [right of=C8] {$k(\sqrt{-1}\,)$};
  \draw (G) -- (C8);
  \draw (G) -- (C81);
  \draw (G) -- (C24);
  \node (C4) [above of=C8] {$k_\gen$};
  \draw (C8) -- (C4);
  \draw (C81) -- (C4);
  \draw (C24) -- (C4);
  \node (H) [above of=C24] {$k(\sqrt{\pi\rho'}\,)$};
  \node (H2) [right of=H]  {$k(\sqrt{\pi\rho}\,)$};
  \draw (C24) -- (H);
  \draw (C24) -- (H2);
  \node (C2) [above of=H] {$k_\gen(\sqrt{\pi\rho}\,)$};
  \node (C21) [above of=H2] {$k(\sqrt{\pi\mu}\,)$};
  \node (C22) [right of=C21] {$k(\sqrt{\rho\mu}\,)$};
  \draw (C4) -- (C2);
  \draw (H) -- (C2);
  \draw (H2) -- (C2);
  \draw (H2) -- (C21);
  \draw (H2) -- (C22);
  \node (E) [above of=C21] {$L$};
  \draw (C2) -- (E);
  \draw (C21) -- (E);
  \draw (C22) -- (E);
\end{tikzpicture} \hskip -2.4cm
\begin{tikzpicture}[node distance=1.5cm]
  \node (G) at (0,0) {$[4,2]$};
  \node (C8)  [above of =G] {$[8]$};
  \node (C81) [left  of=C8] {$[8]$};
  \node (C24) [right of=C8] {$[4,2]$};
  \draw (G) -- (C8);
  \draw (G) -- (C81);
  \draw (G) -- (C24);
  \node (C4) [above of=C8] {$[4]$};
  \draw (C8) -- (C4);
  \draw (C81) -- (C4);
  \draw (C24) -- (C4);
  \node (H) [above of=C24] {$[4]$};
  \node (H2) [right of=H]  {$[2,2]$};
  \draw (C24) -- (H);
  \draw (C24) -- (H2);
  \node (C2) [above of=H] {$[2]$};
  \node (C21) [above of=H2] {$[2]$};
  \node (C22) [right of=C21] {$[2]$};
  \draw (C4) -- (C2);
  \draw (H) -- (C2);
  \draw (H2) -- (C2);
  \draw (H2) -- (C21);
  \draw (H2) -- (C22);
  \node (E) [above of=C21] {$[1]$};
  \draw (C2) -- (E);
  \draw (C21) -- (E);
  \draw (C22) -- (E);
\end{tikzpicture}
\caption{Subfield diagram of the $2$-class field tower of $k$ and
$2$-class groups of its subextensions}
\end{figure}

\begin{table}
  $$ \begin{array}{rr|cc|c}
  \rsp  p &   q & \pi \rho & \pi \rho' & \pi \mu \\ \hline
  \rsp  5 &  13 &  -7 +  4i &   1 +  8i & (2+2i + \sqrt{-7+4i}\,)(1+2i) \\
  \rsp  5 &  53 &  -3 - 16i & -11 + 12i
          & (6+6i + (1+2i)\sqrt{-3-16i}\,)(1+2i) \\   
  \rsp  5 & 197 & -27 - 16i &  29 + 12i
          & (24+24i + (1+8i)\sqrt{-27-16i}\,)(1+2i) \\  
  \rsp 29 &  37 &  17 + 28i & -7 + 32i 
          & (38 + 38i +(11-4i) \sqrt{17+28i}\,)(5+2i) \\
  \rsp  5 & 317 &  17 + 36i & -39 +  8i
          & (12+12i + (9+14i)\sqrt{17+36i}\,)(1+2i) \\
  \rsp  5 & 397 & - 7 + 44i & -31 + 32i
          & (28+28i +(5+2i)\sqrt{-7+44i}\,)(1+2i) \\
  \rsp 37 &  61 & -31 + 36i &  41 + 24i
          & (4+4i + (-1+2i)\sqrt{-31+36i}\,)(1+6i) \\ 
  \rsp 29 & 101 & -15 + 52i &  25 + 48i
          & (80+80i + (-11+6i)\sqrt{-15+52i}\,)(5+2i) \\
  \rsp 13 & 229 &  41 + 36i &  49 + 24i
          & (36+36i + (9-2i)\sqrt{41+36i}\,)(3+2i)
  \end{array} $$
  \caption{Unramified $\MD_4(2)$-extensions of
    $k = \Q(\sqrt{-pq}\,)$.}\label{TabM42}
\end{table}

Observe that $K(\sqrt{-3+16i}\,)$ and $k(\sqrt{-3-16i}\,)$
are isomorphic quadratic extensions of $K = \Q(i,\sqrt{pq}\,)$;
the choice of signs in $\sqrt{-3 \pm 16i}$ in
$$ \mu \rho = (6+6i + (1+2i)\sqrt{-3-16i}\,)(1+2i) $$
is chosen in such a way that we get a solution of (\ref{Epirho}),
so either
$$ (6+6i)^2 - (1+2i)^2 \cdot (-3-16i) = (1+2i)^2 \cdot (-3-16i) $$
or
$$ (6-6i)^2 - (1-2i)^2 \cdot (-3+16i) = (1-2i)^2 \cdot (-3+16i). $$

We remark in passing that combining Kaplan's condition with the
fact that $16 \mid h$ is equivalent to $(\frac2{A+B}) = 1$ implies that 
we have $\Cl_2(k) \simeq [4,2]$ if and only if
$[\frac{\pi}{\rho}] = -(\frac2{A+B})$.

\end{document}